\documentclass{article}[12pt]
\usepackage{amssymb}
\usepackage{epsfig}
\usepackage{amsmath}
\usepackage{amsfonts}
\usepackage{pgf,tikz}

\newcommand{\R}{{{\Bbb R}}}

\newcommand{\N}{{{\Bbb N}}}

\newtheorem{theorem}{\sc Theorem}[section]
\newtheorem{proposition}{\sc Proposition}[section]

\newtheorem{definition}{\sc Definition}[section]
\newtheorem{remark}{\sc Remark}[section]
\newtheorem{corollary}{\sc Corollary}[section]
\newtheorem{example}{\sc Example}[section]

\def\qed{\hbox to 0pt{}\hfill$\rlap{$\sqcap$}\sqcup$\medbreak}

\title{Integration by parts and by substitution unified, Green's Theorem and uniqueness for ODEs}

\author{ Jos\'e \'Angel Cid$^\dag$ and Rodrigo L\'opez Pouso$^\ddag$\footnote{J. \'A. Cid was supported by project EM2014/032, Conseller\'{\i}a de Cultura, Educaci\'on e Ordenaci\'on Universitaria, Xunta de Galicia, Spain. R. L. Pouso  was partially supported by
Ministerio de Econom\'{\i}a y Competitividad, Spain, and FEDER, Project
MTM2013-43014-P. }} 

\date{}
\begin{document}
 \maketitle

\begin{center}  {\small $\dag$ Departamento de Matem\'aticas, Universidade de Vigo, \\ 32004, Pabell\'on 3, Campus de Ourense, Spain.\\ Email: angelcid@uvigo.es}  \end{center}

\begin{center}  {\small $^\ddag$ Departamento de An\'alise Matem\'atica, Universidade de Santiago de Compostela, \\ 15782, Facultade de Matem\'aticas, Campus Vida, Santiago, Spain.\\  Email: rodrigo.lopez@usc.es}
\end{center}

\medbreak

 \abstract{We present a rather unknown version of the change of variables formula for non-autonomous functions. We will show that this formula is equivalent to Green's Theorem for regions of the plane bounded by the graphs of two continuously differentiable functions. Besides, the formula has interesting applications in the uniqueness of solution of ordinary differential equations.}

\medbreak


\medbreak


\section{Introduction}

The Leibniz rule for differentiating an integral, that you can find for instance in the Monthly all-star \cite{flanders}, 
\begin{equation}\label{lr} \frac{d}{dt} \left( \int_{x_1(t)}^{x_2(t)} f(t,r)\, dr \right)=f(t,x_2(t))x_2'(t)-f(t,x_1(t))x_1'(t)+\int_{x_1(t)}^{x_2(t)}{\dfrac{\partial f}{\partial t}(t,r) \, dr} , \end{equation}
is  a gem that includes not only the Fundamental Theorem of Calculus  (FTC)
$$\frac{d}{dt} \left( \int_{a}^{t} f(r)\, dr \right)=f(t),$$
 and the usual differentiation formula under the integral sign 
 $$\frac{d}{dt} \left( \int_{a}^{b} f(t,r)\, dr \right)=\int_{a}^{b}{\dfrac{\partial f}{\partial t}(t,r) \, dr}, $$
  but as well hides the following interesting version of  the change of variables formula:

\begin{theorem}[Non-autonomous change of variables] \label{thmain} Let $f:[a,b]\times [c,d]\subset \R^2\to \R$ be a continuous function with continuous partial derivative with respect to the first variable and let $x:[a,b]\to [c,d]$ be continuously differentiable. Then, 
\begin{equation}\label{change}\int_a^b{f(t,x(t)) \, x'(t) \, dt}=\int_{x(a)}^{x(b)}{f(b,r) \, dr}-\int_a^b {\left(\int_{x(a)}^{x(t)}{\dfrac{\partial f}{\partial t}(t,r) \, dr} \right) dt }. \end{equation}

\end{theorem} 

\noindent {\bf DIY Proof.} Take $x_1(t)\equiv x(a)$ and $x_2(t)=x(t)$ in (\ref{lr}), integrate between $a$ and $b$, apply the FTC and rearrange the terms.  \qed

\begin{remark}
\label{lip}
Formula (\ref{change}) is also valid if $x:[a,b]\to [c,d]$ is Lipschitz continuous or, more generally, absolutely continuous on $[a,b]$. See \cite{ccp}, where the authors proved (\ref{change}) under weaker assumptions and used it to derive existence results for a class of nonlinear differential equations.
\end{remark}

 It is worthy of remark that formula (\ref{change}) unifies the two most important techniques in elementary integration: the usual change of variables formula when $f(t,x)\equiv f(x)$, namely,
 $$\int_a^b{f(x(t)) \, x'(t) \, dt}=\int_{x(a)}^{x(b)}{f(r) \, dr},$$
 and the integration by parts whenever $f(t,x)\equiv f(t)$, that is,
$$  \int_a^b{f(t) \, x'(t) \, dt}=f(b)x(b)-f(a)x(a)-\int_a^b f'(s) x(s) ds .$$
Moreover, as we shall see at Section 2, formula \eqref{change} is equivalent to a version of Green's Theorem for special regions of the plane. Finally, Section 3 contains the main contributions in this paper: we shall show how (\ref{change}) leads to uniqueness criteria for the initial value problem
\begin{equation}
\label{eq1}
x'(t)=f(t,x(t)), \quad x(t_0)=x_0.
\end{equation}

To anticipate the flavor of our uniqueness results we present now a particular case, where we show that a suitable bound on the nonlinearity is enough to get uniqueness of the constant solution.
\begin{theorem}
\label{th3}
Let $f:U=[t_0-a,t_0+a] \times [x_0-b,x_0+b] \longrightarrow \R$ be continuous on $U$.

If there exists $\psi :[0,+\infty)\to[0,+\infty)$ continuous such that 
$\displaystyle\int_{0+} \frac{d\tau}{\psi(\tau)}=+\infty$ and 
$$|f(t,x)| \le   \psi(|x-x_0|)\quad \mbox{for all $(t,x) \in U$,}$$
 then (\ref{eq1}) has a unique solution, namely, $x(t)=x_0$ for all $t \in [t_0-a,t_0+a]$.
\end{theorem} 
Notice that Theorem \ref{th3} extends to non-autonomous equations the sufficient part of the following known result (see \cite{b,w}): if $f:[x_0-b,x_0+b]\to [0,\infty)$ is continuous and $f(x_0)=0$ then problem 
$$x'(t)=f(x(t))\quad t\ge t_0, \quad x(t_0)=x_0,$$
 has a unique  local solution defined on the right of $t_0$ (namely, the constant $x_0$) if and only if $\displaystyle \int_{x_0^+} \displaystyle\frac{ds}{f(s)}=+\infty$.

 \section{An equivalent statement: Green's Theorem}
Formula (\ref{change}) yields a particularly easy proof of the following version of Green's Theorem.
We have used this proof in lectures to sophomores in Mathematics and we think it is particularly appropriate for introducing Green's Theorem at that level.

\begin{theorem}
\label{thgreen}
Let $a,b \in \R$, $a<b$, and let $\varphi, \, \psi \in {\cal C}^1([a,b])$ be such that $\varphi(t) \le  \psi(t)$ for all $t \in [a,b]$.

Consider the plane region $D=\{ (t,x) \in \R^2 \, : \, t \in [a,b], \, \, \varphi(t) \le x \le \psi(t) \}$ and let $\Gamma$ denote its boundary with positive orientation.

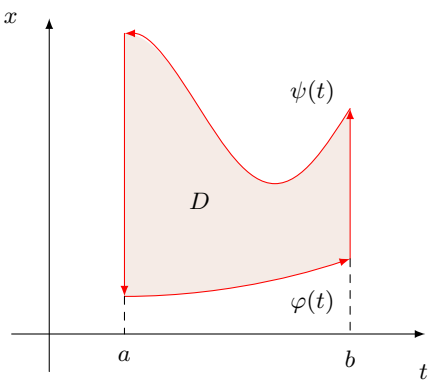
\begin{center}
\definecolor{zzttqq}{rgb}{0.6,0.2,0}
\definecolor{cqcqcq}{rgb}{0.75,0.75,0.75}

\begin{tikzpicture}[line cap=round,line join=round,x=1.0cm,y=1.0cm] 

\draw[-latex,color=black] (-0.5,0) -- (5,0);
\draw (4.8,-0.5) node[anchor=west] {$t$};

\draw[-latex,color=black] (0,-0.5) -- (0,4.2);
\draw (-0.3,4.2) node[anchor=east] {$x$};

\draw[-latex,color=red](1,0.5) parabola (4,1); 
\draw[-latex,color=red]  (4,1) -- (4,3);
\draw[latex-,color=red]  (1,4) cos (2,3) sin (3,2) cos (4,3);
\draw[-latex,color=red]  (1,4) -- (1,0.5);

\fill[color=zzttqq,fill=zzttqq,fill opacity=0.1]   (1,0.5) parabola (4,1) -- (4,3) sin (3,2) cos (2,3) sin (1,4)--(1,0.5);

\draw (3.5,0.7) node[anchor=north] {$\varphi(t)$};
\draw (3.5,3.5) node[anchor=north] {$\psi(t)$};
\draw (2,2) node[anchor=north] {$D$};

\draw[color=black,dashed] (1,0.5) -- (1,0);
\draw (1,-0.1) node[anchor=north] {$a$};
\draw[color=black,dashed] (4,1) -- (4,0);
\draw (4,-0.1) node[anchor=north] {$b$};

\end{tikzpicture}

\end{center}

Assume $D\subset [a,b]\times [c,d]$ and $f_1,f_2: [a,b]\times [c,d]\to \R$ are such that $\dfrac{\partial f_1}{\partial x} $ and $\dfrac{\partial f_2}{\partial t}$ exist and are continuous on $[a,b]\times [c,d]$. Then,  for $F=(f_1,f_2)$, the following identity holds
\begin{equation}
\label{gf}
\int_{\Gamma}F=\int_{D}\left(\dfrac{\partial f_2}{\partial t}-\dfrac{\partial f_1}{\partial x} \right) \, dt dx.
\end{equation}
\end{theorem}

\noindent
{\bf Proof.} By definition  
$$\int_{\Gamma}{F}=\int_{\alpha}^{\beta}{\langle F(\gamma(t)),\gamma'(t) \rangle \, dt},$$
 where $\gamma(\cdot)$ is any positively oriented parametrization of $\Gamma$ and $\langle \cdot,\cdot \rangle$ means scalar product. Hence, it suffices to evaluate and add the following four terms, each one corresponding to a side of $\Gamma$:
 
 \begin{eqnarray}
\label{i1}
\int_a^b{f_1(t,\varphi(t)) \, dt}+\int_a^b{f_2(t,\varphi(t)) \, \varphi'(t) \, dt}, \\
\label{i2}
\int_{\varphi(b)}^{\psi(b)}f_2(b,x) \, dx, \\
\label{i3}
-\int_a^b{f_1(t,\psi(t)) \, dt}-\int_a^b{f_2(t,\psi(t)) \, \psi'(t) \, dt}, \\
\label{i4}
-\int_{\varphi(a)}^{\psi(a)}f_2(a,x) \, dx.
\end{eqnarray}
Using (\ref{change}) in (\ref{i1}) and (\ref{i3}) we get
\begin{align*}
\int_{\Gamma}F &=\int_a^b{f_1(t,\varphi(t)) \, dt}+\int_{\varphi(a)}^{\varphi(b)}{f_2(b,x) \, dx}-\int_a^b \left( \int_{\varphi(a)}^{\varphi(t)}\dfrac{\partial f_2}{\partial t}(t,x) \, dx \right) \, dt \\
& \quad +\int_{\varphi(b)}^{\psi(b)}f_2(b,x) \, dx \\
& \quad -\int_a^b{f_1(t,\psi(t)) \, dt}-\int_{\psi(a)}^{\psi(b)}{f_2(b,x) \, dx}+\int_a^b \left( \int_{\psi(a)}^{\psi(t)}\dfrac{\partial f_2}{\partial t}(t,x) \, dx \right) \, dt \\
& \quad -\int_{\varphi(a)}^{\psi(a)}f_2(a,x) \, dx
\end{align*}
Notice that for the first and fifth terms in this sum the fundamental theorem of calculus yields
$$\int_a^b{f_1(t,\varphi(t)) \, dt}-\int_a^b{f_1(t,\psi(t)) \, dt}=-\int_a^b \left(\int_{\varphi(t)}^{\psi(t)}\dfrac{\partial f_1}{\partial x}(t,x) \, dx \right)dt.$$
Adding the second, fourth, sixth and eighth terms we obtain
\begin{align*}
  \int_{\varphi(a)}^{\varphi(b)}{f_2(b,x) \, dx}  &+ \int_{\varphi(b)}^{\psi(b)}f_2(b,x) \, dx -\int_{\psi(a)}^{\psi(b)}{f_2(b,x) \, dx} -\int_{\varphi(a)}^{\psi(a)}f_2(a,x) \, dx \\
    & \quad = \int_{\varphi(a)}^{\psi(a)}f_2(b,x) \, dx-\int_{\varphi(a)}^{\psi(a)}f_2(a,x) \, dx\\
    & \quad = \int_{\varphi(a)}^{\psi(a)}[f_2(b,x)-f_2(a,x)]  \, dx
\end{align*}
On the other hand, adding the third and seventh terms we get
\begin{align*}
-\int_a^b \left( \int_{\varphi(a)}^{\varphi(t)}\dfrac{\partial f_2}{\partial t}(t,x) \, dx \right) dt  &+
\int_a^b \left( \int_{\psi(a)}^{\psi(t)}\dfrac{\partial f_2}{\partial t}(t,x) \, dx \right) \, dt\\
&=\int_a^b \left( \int_{\varphi(t)}^{\psi(t)}\dfrac{\partial f_2}{\partial t}(t,x) \, dx-
\int_{\varphi(a)}^{\psi(a)}\dfrac{\partial f_2}{\partial t}(t,x) \, dx \right) dt \\
&=\int_a^b \left( \int_{\varphi(t)}^{\psi(t)}\dfrac{\partial f_2}{\partial t}(t,x) \, dx \right) dt-\int_{\varphi(a)}^{\psi(a)}\int_a^b \dfrac{\partial f_2}{\partial t}(t,x) \, dt \, dx \\
&=\int_a^b \left( \int_{\varphi(t)}^{\psi(t)}\dfrac{\partial f_2}{\partial t}(t,x) \, dx \right) dt -\int_{\varphi(a)}^{\psi(a)}[f_2(b,x)-f_2(a,x)] \, dx 
\end{align*}
Finally, adding all the terms we obtain the desired result
$$
\int_{\Gamma}F=\int_a^b \left( \int_{\varphi(t)}^{\psi(t)}\dfrac{\partial f_2}{\partial t}(t,x) \, dx \right) dt -\int_a^b \left(\int_{\varphi(t)}^{\psi(t)}\dfrac{\partial f_1}{\partial x}(t,x) \, dx \right)dt=\int_{D}\left(\dfrac{\partial f_2}{\partial t}-\dfrac{\partial f_1}{\partial x} \right) \, dt dx.
$$

\qed

In fact Theorem \ref{thgreen} implies Theorem \ref{thmain}, and thus both are equivalent. Indeed, consider $f(t,x)$ and $x(t)$ as in the assumptions of Theorem \ref{thmain} and define $\varphi(t)\equiv \displaystyle\min_{t\in [a,b]} x(t)$, $\psi(t)=x(t)$, $f_1\equiv 0$ and $f_2=f$. Let us call $c=\displaystyle\min_{t\in [a,b]} x(t)$, then we have
\begin{equation}\label{eqg1}  \int_{\Gamma}F =   \int _c^{x(b)} f(b,r) dr-\int_a^b f(t,x(t))x'(t)dt-\int_c^{x(a)}f(a,r)dr 
 \end{equation}
 and
 
 \begin{align}
 \label{eqg2}  
\int_{D}\left(\dfrac{\partial f_2}{\partial t}-\dfrac{\partial f_1}{\partial x} \right) \, dt dx & = \int_a^b\int_c^{x(t)} \dfrac{\partial f} {\partial t}(t,r) \, dr dt \notag \\
    & =  \int_a^b\int_c^{x(a)} \dfrac{\partial f} {\partial t}(t,r) \, dr dt+ \int_a^b\int_{x(a)}^{x(t)} \dfrac{\partial f} {\partial t}(t,r) \, dr dt  \notag\\  
    & = \int_c^{x(a)}  \int_a^b \dfrac{\partial f} {\partial t}(t,r) \, dt dr+ \int_a^b\int_{x(a)}^{x(t)} \dfrac{\partial f} {\partial t}(t,r) \, dr dt\notag \\  
   & = \int_c^{x(a)} [f(b,r)-f(a,r)]\, dr+ \int_a^b\int_{x(a)}^{x(t)} \dfrac{\partial f} {\partial t}(t,r) \, dr dt.
\end{align}
 
Since \eqref{eqg1} and \eqref{eqg2} are equal by Theorem \ref{thgreen}, we obtain 
$$\int_a^b f(t,x(t))x'(t)dt=\int _{x(a)}^{x(b)} f(b,r) dr- \int_a^b\int_{x(a)}^{x(t)} \dfrac{\partial f} {\partial t}(t,r) \, dr dt,$$
which is formula  (\ref{change}), as desired.

\section{Uniqueness criteria for ODE's}

Uniqueness for differential equations is an old subject far from being solved, see \cite{agalak,hartman} and references therein, and that still sparks interest in researching \cite{cp3,cons,dns,ferreira,h}. In this section we show how our general version of the formula of change of variables yields new effective conditions for uniqueness.

Our approach to uniqueness recaptures original ideas by Perron and Kamke, and it is based on the following definition. Here and henceforth, $a$ and $b$ are positive real numbers.

\begin{definition}
\label{defad}
A function $g:(0,a] \times [0,b] \longrightarrow [0,+\infty)$ is a uniqueness bound if for any $\tilde a \in (0,a]$ the unique Lipschitz continuous function $\varphi:[0,\tilde a] \longrightarrow [0,+\infty)$ such that 
\begin{equation}
\label{subcd}
\varphi'(t) \le g(t,\varphi(t)) \, \, \mbox{for almost all (a.a.) $t \in (0,\tilde a]$, and}\quad \varphi(0)=\varphi'(0)=0,
\end{equation}
is $\varphi(t)=0$ for all $t \in [0,\tilde a]$.
 \end{definition}

Obviously, $g(t,x)=0$ for all $(t,x)$ is a uniqueness bound, and we shall discover many nontrivial ones soon by means of our general formula of change of variables. Notice that uniqueness bounds need not be continuous functions. 

\subsection{What are uniqueness bounds good for?}

Uniqueness bounds can be used to prove uniqueness of solutions for both ordinary differential equations and systems of ODEs as we show in our next proposition. We remark that arguments are simpler for equations: in particular,  we can use a simpler version of Definition \ref{defad} which uses continuously differentiable functions instead of Lipschitz continuous ones, and replaces (\ref{subcd}) by
\begin{equation}
\label{subcd2}
\varphi'(t) \le g(t,\varphi(t)) \, \, \mbox{for all $t \in (0,\tilde a]$, and}\quad \varphi(0)=\varphi'(0)=0.
\end{equation}
We shall indicate how we can simplify proofs for equations at relevant places.

From now on, we assume that $t_0 \in \R$ and $x_0=(x_{0,1},x_{0,2},\dots,x_{0,n}) \in \R^n$, $n \in \N$, are fixed, $\| \cdot \|$ denotes a norm in $\R^n$ (the specific form of which is not important), and we discuss uniqueness of solutions for the initial value problem
\begin{equation}
\label{ivpub1}
x'=f(t,x), \quad x(t_0)=x_0.
\end{equation}
We shall denote by $\overline{B(x_0,r)}$ the closed ball centered at $x_0$ and radius $r>0$ corresponding to the previous norm.

The second part of the following proposition is a variant of Kamke's uniqueness theorem which is enough for our objectives. With the aid of uniqueness bounds we have been able to avoid in its proof the use of differential inequalities, providing in this way a self-contained proof  considerably simpler than those available in the literature, cf. \cite[Theorem 2.3]{cl} or \cite[Theorem 6.1]{hartman}. 
\begin{proposition}
\label{proub}
Let $g:(0,a] \times [0,b] \longrightarrow [0,+\infty)$ be a uniqueness bound. Then the following statements are true:
\begin{enumerate}
\item If $f:V = [t_0-a,t_0+a]\times \overline{B(x_0,b)} \longrightarrow \R^n$ is continuous and  
\begin{equation}
\label{ub1}
\|f(t,x)\| \le g( |t-t_0|,\|x-x_0\|) \quad \mbox{for all $(t,x) \in V$, $t\neq t_0$,}
\end{equation}
and $f(t,x_0)=(0,0,\dots,0)$ for all $t \in [t_0-a,t_0+a]$,
then the initial value problem (\ref{ivpub1}) has only the constant solution $x \equiv x_0$ on any subinterval $I \subset [t_0-a,t_0+a]$.
\item If $f:W=[t_0-a,t_0+a]\times \overline{B(x_0,b/2)}\longrightarrow \R^n$ is continuous and
\begin{equation}
\label{ub2}
\|f(t,x)-f(t,y)\| \le g( |t-t_0|,\|x-y\|) \quad \mbox{for all $(t,x), (t,y) \in W$, $t\neq t_0$,}
\end{equation}
then the initial value problem (\ref{ivpub1}) has at most one solution on any subinterval $I \subset [t_0-a,t_0+a]$.
\end{enumerate}
\end{proposition}

\noindent
{\bf Proof.} We are going to prove both claims at one stroke: assume, reasoning by contradiction, that we can find two different solutions of (\ref{ivpub1}), say $x,\, y:I \subset [t_0-a,t_0+a] \longrightarrow \R^n$ (w.l.g. we assume that $y(t)=x_0$ on $I$ in the first case). Then we can find $t_1, t_2 \in (-a,a)$ such that $x(t_0+t_1)=y(t_0+t_1)$ and one of the following situations holds:
\begin{enumerate}
\item[$(a)$] $0 \le t_1 < t_2$ and $x(t+t_0) \neq y(t+t_0)$ for all $t \in (t_1,t_2]$; or
\item[$(b)$] $t_2<t_1 \le 0$ and $x(t+t_0) \neq y(t+t_0)$ for all $t \in [t_2,t_1)$.
\end{enumerate}
Let us prove that $(a)$ leads to a contradiction. We define a real valued function as follows: $\varphi (t)=0$ for all $t \in [0,t_1]$ and $\varphi(t)=\|x(t+t_0)-y(t+t_0)\|$ for all $t \in [t_1,t_2]$. We are going to show that it satisfies
all the conditions in Definition \ref{defad}\footnote{In dimension $n=1$ we can prove that $\varphi \in {\cal C}^1([0,t_2])$ and satisfies (\ref{subcd2}) on $(0,\tilde a]=(0,t_2]$.}. First, note that for all $s,t \in [t_1,t_2]$, $s<t$, we have
\begin{align}
\nonumber
|\varphi(t)-\varphi(s)|& =| \, \|x(t+t_0)-y(t+t_0)\|-\|x(s+t_0)-y(s+t_0)\| \, | \\
\nonumber
& \le \| x(t+t_0)-x(s+t_0)+y(s+t_0)-y(t+t_0)\| \\
\nonumber
&= \left\| \int_{s+t_0}^{t+t_0}[x'(r)-y'(r)] \, dr \right\| =\left\| \int_{s+t_0}^{t+t_0}[f(r,x(r))-f(r,y(r))]  \, dr \right\| \\
\label{iinn}
& \le (t-s) \, \max_{s+t_0 \le r \le t+t_0}\|f(r,x(r))-f(r,y(r))\|,
\end{align}
which implies that $\varphi$ is Lipschitz continuous on $[t_1,t_2]$ with Lipschitz constant
$$L=\max_{t_1+t_0 \le r \le t_2+t_0}\|f(r,x(r)) -f(r,y(r))\|.$$
Furthermore, $\varphi$ is Lipschitz continuous on the whole of $[0,t_2]$ because it is constant on $[0,t_1]$ and continuous at $t_1$.
In particular, $\varphi'(s)$ exists for almost all $s \in [0,t_2]$, and for any
of those points $s \in [t_1,t_2)$, we deduce from (\ref{iinn}) and either condition (\ref{ub1}) or (\ref{ub2}) depending on the case (remember that in the first case $y(t)\equiv x_0$) that
$$\varphi'(s) \le |�\varphi'(s)| = \lim_{t \to s^+} \left| \dfrac{\varphi(t)-\varphi(s)}{t-s} \right| \le \|f(s+t_0,x(s+t_0))-f(s+t_0,y(s+t_0))\| \le g(s,\varphi(s)).$$
Therefore, $\varphi'(s) \le g(s,\varphi(s))$ for almost all $s \in [0,t_2]$.
 
Finally, $\varphi'(0)=0$ if $0 < t_1$, and if $t_1=0$ we deduce the same result by means of (\ref{iinn}) with $s=0$. Indeed, for every $t>0$ we have 
$$ \dfrac{|\varphi(t)|}{t} \le \max_{t_0 \le r \le t+t_0}\|f(r,x(r))-f(r,y(r)) \|,$$
and the assumptions ensure that
$$\lim_{t \to 0^+}\max_{t_0 \le r \le t+t_0}\|f(r,x(r))-f(r,y(r)) \|=\|f(t_0,x(t_0))-f(t_0,y(t_0))\|=0.$$
Hence $\varphi'(0)=0=\varphi(0)$ and $\varphi$ satisfies (\ref{subcd}). This implies that $\varphi \equiv 0$ on $[0,t_2]$ because $g$ is a uniqueness bound, a contradiction with $(a)$.

We obtain a similar contradiction if we assume $(b)$ and we consider $\varphi(t)=\|x(t_0-t)-y(t_0-t)\|$ for $t \in [-t_1,-t_2]$, $\varphi \equiv 0$ on $[0,-t_1]$.  \qed

 \subsection{When is $g(t,x)$ a uniqueness bound?}
 
Next we state and prove the main result in this section, where we show the applicability of the general formula of change of variables for checking that a given function $g(t,x)$ is a uniqueness bound.  
\begin{theorem}
\label{th33}
A function $g:(0,a] \times [0,b] \longrightarrow [0,+\infty)$ is a uniqueness bound provided that $g(t,x) \le p(t) \psi(t,x)$, where $p:(0,a] \longrightarrow (0,+\infty)$ is locally integrable on $(0,a]$, $\psi :(0,a] \times [0,2b]\to[0,+\infty)$ is continuous, and the following properties are satisfied:
\begin{enumerate}
\item[(i)] $\psi(t,0)=0$ for all $t \in (0,a]$ and $\psi(t,x)>0$ for all $(t,x) \in (0,a]\times (0,b]$;
\item[(ii)] $1/\psi$ has a continuous partial derivative with respect to $t$  on $(0,a)\times (0,b)$;  and
\item[(iii)] For every Lipschitz continuous\footnote{In dimension $n=1$ we can restrict our attention to continuously differentiable functions $u:[t_1,t_2] \longrightarrow [0,b]$ which are positive on $(t_1,t_2]$.} function $u:(t_1,t_2] \subset (0,a] \longrightarrow (0,b)$ such that  
\begin{equation}
\label{limit}
\lim_{t \to t_1^+}\dfrac{u(t)}{t-t_1}=0,
\end{equation}
we have
\begin{equation}
\label{inf1}
\limsup_{t \to t_1^+}  \left( \displaystyle\int_{u(t)}^{u(t_2)} \frac{1}{\psi(t_2,r)}dr - \int_{t}^{t_2}{p(s) \, ds}  \\
-\int_{t}^{t_2}{\int_{u(t)}^{u(s)}\dfrac{\partial}{\partial s}\dfrac{1}{\psi(s,r)} \, dr\, ds}\right)>0.
\end{equation}
\end{enumerate}
The result holds valid if we replace conditions (ii) and (iii) by, respectively, 
\begin{enumerate}
 \item[(ii')] $\psi$ is nondecreasing with respect to its first variable; and
\item[(iii')] For every Lipschitz continuous function $u:(t_1,t_2] \subset (0,a] \longrightarrow (0,b)$ satisfying (\ref{limit})  
we have
\begin{equation}
\label{inf2}
\limsup_{t \to t_1^+}  \left( \displaystyle\int_{u(t)}^{u(t_2)} \frac{1}{\psi(t_2,r)}dr - \int_{t}^{t_2}{p(s) \, ds}  \right)>0.
\end{equation}
\end{enumerate} 
 
\end{theorem} 

\noindent
{\bf Proof.} Let $\tilde a�\in (0,a]$ be fixed and let $\varphi:[0,\tilde a] \longrightarrow [0,+\infty)$ be a Lipschitz continuous function satisfying (\ref{subcd}); we have to prove that $\varphi \equiv 0$ on $[0, \tilde a]$.

Reasoning by contradiction, we assume that we can find $t_1 \in [0,\tilde a)$ and $t_2 \in (t_1,�\tilde a]$ such that $\varphi(t)>0$ for all $t \in (t_1,t_2]$ and $\varphi(t_1)=0$. 

Notice that $\varphi$ satisfies (\ref{limit}). Indeed, (\ref{limit}) follows from the second part of (\ref{subcd}) if $t_1=0$. On the other hand, if $t_1>0$ then the Fundamental Theorem of Calculus for the Lebesgue integral guarantess that for all $t \in (t_1,t_2)$  we have
\begin{align*}
0 \le \dfrac{\varphi(t)}{t-t_1} &= \dfrac{1}{t-t_1}\int_{t_1}^t{\varphi'(s) \, ds} \\
& \le  \dfrac{1}{t-t_1}\int_{t_1}^t{p(s) \psi(s,\varphi(s)) \, ds} \\
& \le \dfrac{\max_{t_1 \le s \le t}\psi(s,\varphi(s))}{t-t_1}\int_{t_1}^{t_2}{p(s) \, ds},
\end{align*}
and since $p$ is integrable on $[t_1,t_2]$, $\psi(s,\varphi(s))$ is continuous at $t=t_1$ and $\psi(t_1,\varphi(t_1))=0$, we deduce that $\varphi$ satisfies (\ref{limit}) also in case $t_1>0$.

For almost all $t \in (t_1,t_2]$ we have
\begin{equation}
\label{d}
\varphi'(t) \le g(t,\varphi(t)) \le p(t)\psi(t,\varphi(t)),
\end{equation} so the formula of change of variables \eqref{change} (see Remark \ref{lip}) yields
 \begin{align*}
\int_{t}^{t_2}{p(s) \, ds}&\ge \int_{t}^{t_2}{\dfrac{\varphi'(s)}{\psi(s,\varphi(s))}ds}\\
&= \int_{\varphi(t)}^{\varphi(t_2)}{\dfrac{dr}{\psi(t_2,r)}} \\
& \mbox{} \quad -\int_t^{t_2} \int_{\varphi(t)}^{\varphi(s)}\dfrac{\partial}{\partial s}\dfrac{1}{\psi(s,r)} \, dr \,ds.\\
\end{align*}
Taking limits when $t \to t_1^+$ on the previous inequality leads to a contradiction with (\ref{inf1}) with $u(t)=\varphi(t)$.  

To prove that the result holds valid when we replace $(ii)$ and $(iii)$ by $(ii')$ and $(iii')$ simply repeat the proof and notice that we can infer from (\ref{d}) and $(ii')$ that
$$\varphi'(t) \le p(t)\psi(t_2,\varphi(t)) \quad \mbox {for a.a. $t \in (t_1,t_2),$}$$
and then for every fixed $t \in (t_1,t_2)$ we have
$$
\int_{t}^{t_2}{p(s) \, ds}  \ge \int_{t}^{t_2}{\dfrac{\varphi'(s)}{\psi(t_2,\varphi(s))}ds} = \int_{\varphi(t)}^{\varphi(t_2)}{\dfrac{dr}{\psi(t_2,r)}},$$
where we have used the usual formula of change of variables. Now the proof follows exactly as in the previous case.
 \qed

The use of different $p$'s and $\psi$'s in Theorem \ref{th33} yields different uniqueness results via the second part of Proposition \ref{proub}. Some examples follow easily (check conditions (ii') and (iii') in all of them):
\begin{enumerate}
\item (Lipschitz's Theorem) $p=1$ and $\psi(t,x)=c \, x$, where $c>0$ is fixed.
\item (Osgood's Theorem) $p=1$ and $\psi(t,x)=\psi(x)$, where $\psi(x)>0$ for $x>0$ and $\displaystyle\int_{0^+}{\frac{ds}{\psi(s)}}=+\infty$.
\item (Montel--Tonelli's Theorem)  $\displaystyle\int_{0^+}{p(s) \, ds}<+\infty$, and $\psi(t,x)=\psi(x)$ as in Osgood's Theorem.
\item (Nagumo's Theorem) $p(t)=1/t$ and $\psi(t,x)= x$.
\item (Van Kampen's Theorem, \cite{vK})  $p(t)=(1+q(t))/t$, where $q(t)\ge 0$ for $t>0$ and $\displaystyle\int_{0}^{a}{\frac{q(s)}{s}ds}<+\infty$,  and $\psi(t,x)=x$.
\end{enumerate}

Notice that every uniqueness result in the previous list can be proven by means of suitable separable uniqueness bounds, that is of the form $g(t,x)=p(t)\psi(x)$, and this was essentially known by LaSalle \cite{agalak,lasalle}. Indeed, following the presentation in \cite[Corollary 1.15.6]{agalak} the second part of Proposition \ref{proub} holds true with $g(t,x)=p(t)\psi(x)$ provided that $p(t)\ge 0$ is continuous for $t>0$, $\psi(x)$ is continuous for $x\ge 0$, $\psi(0)=0$, $\psi(x)>0$ for $x>0$ and either one of the following conditions holds:

\begin{equation}\label{admis1} \limsup_{t\to 0^+} \int_t \left( \frac{1}{\psi(s)} -p(s)\right)ds=\infty; \quad \mbox{or}
\end{equation}

\begin{equation}\label{admis2} \limsup_{t\to 0^+} \int_t \left( \frac{1}{\psi(s)} -p(s)\right)ds>-\infty \quad \mbox{and} \quad \psi(x)\le x.
\end{equation}
\medbreak \medbreak

It is easy to check that each one of LaSalle's conditions \eqref{admis1} or \eqref{admis2} implies condition {\it (iii')} in Theorem \ref{th33}. On the other hand, the functions $p(t)=\displaystyle\frac{1}{t}$ and $\psi(x)={\rm e}^x-1$ do not satisfy either  \eqref{admis1} or \eqref{admis2}, but they can be bounded above by functions satisfying (\ref{inf1}), as we
shall show as a consequence of the following result, whose proof leans on our general formula \eqref{change} in its full form.  

\begin{corollary}
\label{co1}
A function $g:(0,a] \times [0,b] \longrightarrow [0,+\infty)$ is a uniqueness bound provided that $g(t,x) \le p(t) \psi(t,x)$  for $$p(t)=(1+q_1(t))/t \quad \mbox{and} \quad \psi(t,x)=(1+q_2(t)x^{\gamma})x,$$ where $q_1:(0,a] \longrightarrow [0,+\infty)$ is measurable, $\displaystyle\int_{0}^{a}{\frac{q_1(s)}{s}ds}<+\infty$, $q_2:[0,a] \longrightarrow [0,+\infty)$ is continuous, $q_2'$ is continuous and integrable on $(0,a)$, and $\gamma >0$.
\end{corollary}

\noindent
{\bf Proof.}  Obviously, conditions $(i)$ and $(ii)$ in Theorem \ref{th33} are satisfied. Let us check condition $(iii)$: we consider a Lipschitz continuous function $u:(t_1,t_2]\subset (0,a] \longrightarrow (0,b)$ satisfying (\ref{limit}) and we have to show that it satisfies (\ref{inf1}).

For each $t \in (t_1,t_2)$ we compute
\begin{align*}
\displaystyle\int_{u(t)}^{u(t_2)} \frac{dr}{\psi(t_2,r)}&-\int_{t}^{t_2}{p(s) \, ds}-   
\int_{t}^{t_2}{\int_{u(t)}^{u(s)}\dfrac{\partial}{\partial s}\dfrac{1}{\psi(s,r)} \, dr\, ds} \\
&=\int_{u(t)}^{u(t_2)}{\dfrac{dr}{(1+q_2(t_2)r^{\gamma})r}} \\
& \quad -\mbox{ln}\dfrac{t_2}{t}-\int_{t}^{t_2}{\dfrac{q_1(s)}{s}ds} \\
& \quad +\int_{t}^{t_2}{\int_{u(t)}^{u(s)}\dfrac{{q_2}'(s)r^{\gamma-1}}{(1+q_2(s)r^{\gamma})^2} \, dr\, ds}
\end{align*}
The assumptions guarantee that the third and fourth terms are bounded as functions of $t \in (t_1,t_2)$. For the remaining two terms we have
$$
\int_{u(t)}^{u(t_2)}{\dfrac{dr}{(1+q_2(t_2)r^{\gamma})r}}  -\mbox{ln}\dfrac{t_2}{t} 
= -\mbox{ln}\left(\dfrac{t_2}{t}\dfrac{(1+q_2(t_2)u^{\gamma}(t_2))^{1/\gamma}}{u(t_2)}\dfrac{u(t)}{(1+q_2(t_2)u^{\gamma}(t))^{1/\gamma}} \right),
$$
and this function tends to $+\infty$ as $t$ tends to $t_1$ from the right thanks to (\ref{limit}). Therefore in this case the limit in (\ref{inf1}) is $+\infty$ and the proof is complete. \qed

A consequence of Corollary \ref{co1} improves on Nagumo's uniqueness bound $g(t,x)=x/t$. We emphasize that $g(t,x)=cx/t$ is not a uniqueness bound if $c>1$, since in this case $\varphi(t)=t^{c}$ is a non-trivial Lipschitz continuous solution of $\varphi'(t)=g(t,\varphi(t))$ for $t>0$, $\varphi(0)=\varphi'(0)=0$. Even more, in fact there exist examples of non-uniqueness for the initial value problem \eqref{ivpub1} with $f(t,x)$ satisfying inequality \eqref{ub2} with $g(t,x)=cx/t$, $c>1$ (see \cite[Example 1.6.1]{agalak}). 

\begin{corollary}
\label{co2}
A function $g:(0,a] \times [0,b] \longrightarrow [0,+\infty)$ is a uniqueness bound provided that 
$$g(t,x) \le \dfrac{\varphi(x)}{t} \quad  \mbox{for all $(t,x) \in (0,a]\times [0,b]$,}$$
where $\varphi \in {\cal C}^{1}([0,b])$, $\varphi(0)=0$, $\varphi(x) >0$ for all $x \in (0,b]$, $\varphi'(0) \le 1$, and $\varphi''$ exists and is bounded above on $(0,b)$.
\end{corollary}

\noindent
{\bf Proof.} By Taylor's Theorem, for each fixed $x \in (0,b]$ there exists some $y \in (0,x)$ such that
$$\varphi(x)=\varphi'(0)x+\dfrac{\varphi''(y)}{2}x^2,$$
and therefore the assumptions guaratee the existence of some constant $c>0$ such that
$$\varphi(x) \le x+cx^2 \quad \mbox{for all $x \in [0,b]$.}$$
Hence, $g(t,x) \le p(t)\psi(t,x)$ for $p(t)=1/t$ and $\psi(t,x)=(1+cx)x$, and Corollary \ref{co1} applies. \qed

\begin{example}
Let $a,b  \in (0,+\infty)$ be fixed and consider the function $f:U=[0,a]\times [-b,b]\longrightarrow \R$ defined as
$$f(t,x)=\left\{
\begin{array}{cl} 
  t, & \mbox{if $|x| > \mbox{\rm ln} (1+t^2)$,} \\
  \dfrac{e^{|x|}-1}{t}, & \mbox{if $|x| \le \mbox{\rm ln} (1+t^2)$, $t>0$,}\\
0, &\mbox{if $(t,x)=(0,0)$.}

\end{array}
\right.$$

Obviously, $x \equiv 0$ solves
\begin{equation}
\label{ivpex}
x'=f(t,x), \, \, t \in [0,a], \, \, x(0)=0,
\end{equation}
and we are going to study whether it is the unique solution.

First, notice that \begin{equation}
\label{dd}
0 \le f(t,x) \le t \quad \mbox{for all $(t,x) \in U$,}
\end{equation} 
and therefore $f$ is continuous on $U$. Moreover, for all $(t,x) \in U$, $t>0$, we have
\begin{equation}
\label{cana}
| f(t,x) |=f(t,x) \le \dfrac{\varphi(|x|)}{t},
\end{equation}
for $\varphi(x)=e^x-1$, $x \in [0,b]$, which satisfies the conditions in Corollary \ref{co2}. Hence problem (\ref{ivpex}) has only the zero solution.

Notice that $\varphi(x)>x$ for $x>0$, and therefore we cannot deduce from (\ref{cana}) that a Nagumo condition is satisfied. Morevorer, it is important to note that $f$ does not satisfy any local Lipschitz condition with respect to $x$ or with respect to $t$, and therefore (\ref{ivpex}) falls outside the scope of recent uniqueness results such as those in \cite{dns, h}.
 
\end{example}

\subsection{What about differential inequalities?}

We have already mentioned that with the aid of uniqueness bounds we have been able to avoid the use of differential inequalities in the proof of our version of Kamke's uniqueness theorem. Actually, uniqueness bounds are deeply related with differential inequalities and as an instance we shall present an alternative proof of a particularly famous one: Gronwall's Lemma. For the sake of a clearer presentation, we start with the following elementary observation which gives us another characterization of uniqueness bounds.

\begin{proposition}
\label{proubb}
Let $g:(0,a] \times [0,b]\longrightarrow [0,+\infty)$ be a given function. The following statements are equivalent:
\begin{enumerate}
\item The function $g$ is a uniqueness bound;
\item If $\varphi:[0,\tilde a]\subset[0,a] \longrightarrow \R$ is Lipschitz continuous and
\begin{equation}
\label{ub22}
\varphi'(t) \le g(t,|\varphi(t)|) \quad \mbox{for a.a. $t \in [0,\tilde a]$,} \quad \varphi(0) \le 0, \, \, \varphi'(0) \le 0,
\end{equation}
then $\varphi(t) \le 0$ for all $t \in [0,\tilde a]$.

\end{enumerate}

\end{proposition}

\noindent
{\bf Proof.} We only have to prove that 1 implies 2 because the converse is trivial. Assume that $g(t,x)$ is a uniqueness bound and let $\varphi:[0,\tilde a]\subset[0,a] \longrightarrow \R$ be a Lipschitz continuous
function satisfying (\ref{ub22}). We have to prove that $\varphi \le 0$ on $[0,\tilde a]$. Reasoning by contradiction, we assume that we can find $t_1, t_2
\in [0,\tilde a]$, $t_1 <t_2$, such that 
\begin{equation}
\label{dem1}
\varphi(t) >0  \quad \mbox{for all $t \in (t_1,t_2].$}
\end{equation}
Now we define a function $\phi:[0,t_2]\longrightarrow [0,+\infty)$ as follows: $\phi(t)=0$ for all $t \in [0,t_1]$, and $\phi(t)=\varphi(t)$ for all $t \in (t_1,t_2]$. Obviously, $\phi$ is a nonnegative Lipschitz continuous function, $\phi(0)=0$, and, by (\ref{ub22}) we have
$$\phi'(t) \le g(t,\phi(t)) \quad  \mbox{for a.a. $t \in [0,t_2]$.}$$
Finally, we shall prove that $\phi'(0)=0$. Indeed, the result is trivial if $t_1>0$, and if $t_1=0$ we note that
$$\phi'_+(0)=\varphi'(0) \le 0,$$
and, on the other hand, $\varphi'(0) \ge 0$ thanks to (\ref{dem1}) with $t_1=0$. Hence $\phi'(0)=0$. 

Since $g$ is an uniqueness bound, we deduce that $\phi \equiv 0$ on $[0,\tilde a]$, a contradiction with (\ref{dem1}).
\qed

 Gronwall's Lemma is a particular case of a more general statement about differential inequalities, see \cite[Theorem 6.1]{hale}, but it is interesting by itself. As announced we present a new proof based on uniqueness bounds.

\begin{theorem}
Let $t_0 \in \R$ and $a>0$ be fixed and denote $I=[t_0,t_0+a]$.

If $\alpha, \beta \in {\cal C}(I)$, $\varphi \in {\cal C}^1(I)$, and
\begin{equation}
\label{gl1}
\varphi'(t) \le \alpha(t)+\beta(t)\varphi(t) \quad \mbox{for all $t \in (t_0,t_0+a]$,}
\end{equation}
then 
\begin{equation}
\label{gl2}
\varphi(t) \le \varphi(t_0)\mbox{exp}\left(\int_{t_0}^t\beta(r) \, dr \right)+\int_{t_0}^{t}{\alpha(s) \mbox{exp}\left(\int_s^t{\beta(r) \, dr} \right) ds}
 \quad \mbox{for all $t \in I$.}
\end{equation}
\end{theorem}
\noindent
{\bf Proof.} The right--hand side in (\ref{gl2}), namely
$$x(t)=\varphi(t_0)\mbox{exp}\left(\int_{t_0}^t\beta(r) \, dr \right)+\int_{t_0}^{t}{\alpha(s) \mbox{exp}\left(\int_s^t{\beta(r) \, dr} \right) ds},
\quad t  \in I,$$
solves the linear problem
$$x'=\alpha(t)+\beta(t)x, \,\, t \in I, \, \, x(t_0)=\varphi(t_0).$$

We have to show that the function $\phi(t)=\varphi(t+t_0)-x(t+t_0)$ is nonpositive on $J=[0,a]$. To do it, we note that $\phi$ is Lipschitz continuous on $[0,a]$, $\phi(0)=0$, $\phi'(0)=\varphi'(t_0)-x'(t_0) \le \beta(t_0)(\varphi(t_0)-x(t_0))=0$, and for all $t \in (0,a]$ we have
$$\phi'(t) =\varphi'(t+t_0)-x'(t+t_0) \le \beta(t+t_0)\phi(t) \le g(t,|\phi(t)|),$$
where $g(t,x)=|\beta(t+t_0)|x$.

Theorem \ref{th33} guarantees that $g$ is a uniqueness bound, and then Proposition \ref{proubb} ensures that $\phi \le 0$. The proof is complete.
\qed

\end{document}